\documentclass[a4paper,11pt]{amsart}

\usepackage[francais]{babel}
\usepackage{amssymb} 
\usepackage[latin1]{inputenc}
\usepackage[matrix,arrow]{xy}
\usepackage[mathscr]{eucal}
\usepackage{epic,eepic}
\usepackage{bbm}

\newtheorem{prop}[subsubsection]{Proposition}

\newtheorem{theoreme}[subsubsection]{Théorème}

\newtheorem{lemme}[subsubsection]{Lemme}

\newtheorem{surtheoreme}{Théorème}
\newtheorem{surcorollaire}{Corollaire}
\newtheorem{surprop}{Proposition}

\newtheorem{corollaire}[subsubsection]{Corollaire}

\theoremstyle{definition}

\newenvironment{pf}
        {\medskip\noindent {\it Démonstration --- \ }}
        {\hfill\nobreak $\Box$ \par\bigbreak}

\newcommand{\Hom}{\text{Hom}}

\newcommand{\F}{ \mathbb F  } 
\newcommand{\C}{{ \mathbb C  }}

\newcommand{\Q}{{ \mathbb Q } }
\newcommand{\Z}{{ \mathbb Z  }}

\newcommand{\Ker}{{\text{Ker}\,}}

\renewcommand{\ker}{{\text{Ker}\,}}
\newcommand{\Gal}{{\text{Gal}\,}}
\newcommand{\Disc}{{\text{Disc}\,}}

\newcommand{\anneau}{{ \mathcal O}}

\newcommand{\Id}{{\text{Id}}}
\newcommand{\sss}{{\text{ss}}}

\newcommand{\U}{{\text{U}}}

\newcommand{\Gl}{{\text {GL}}}

\newcommand{\Sl}{{\text {SL}}}

\newcommand{\TT}{{\mathfrak T}}

\newcommand{\rec}{{\text{rec}}}

\newcommand{\St}{{St}}
\newcommand{\st}{{St}}
 
\newcommand{\Hecke}{{\mathcal{H}}}
\newcommand{\ZZ}{{\mathcal{Z}}}

\newcommand{\ad}{{\text{ad}}}
\newcommand{\ssc}{{\text{sc}}}

\newcommand{\Res}{{\text{Res}}}
\newcommand{\ses}{{\text{ss}}}

\newcommand{\pig}{{\varpi}}

\newcommand{\Frob}{{\text{Frob\,}}}
\newcommand{\frob}{{\text{Frob\,}}}

\newcommand{\tr}{{\text{tr\,}}}

\newcommand{\sd}{\rtimes}

\newcommand{\res}{{\text{res}}}

\begin{document}

\baselineskip 14pt

\bibliographystyle{style} 

\title[Sur la compatibilité]{Sur la compatibilité entre les correspondances de Langlands locale et 
globale pour $\U(3)$}
\author{Joël Bellaïche}

\begin{abstract} En utilisant un argument d'augmentation du niveau (et un 
résultat de Larsen sur l'image des représentations galoisiennes 
apparaissant dans des systèmes compatibles), nous prouvons que pour toute 
forme automorphe $\pi$ de $\U(3)$, la représentation galoisienne $l$-adique $\rho_l$ attachée à $\pi$ par Blasius et Rogawski est à chaque place finie
celle associée à $\pi$ par la correspondance de Langlands locale (au
moins à semi-simplification près, et pour un ensemble de densité $1$ de $l$).
Nous nous appuyons sur le travail de Harris et Taylor, qui ont prouvé les mêmes
résultats (pour $\U(n)$) sous l'hypothèse que le changement de base de
$\pi$ est de carré intégrable en au moins une place. Un corollaire de notre
résultat est que toute représentation automorphe pour $G$ tempérée 
à une infinité de places est tempérée partout.  
   
\noindent {\sc Abstract.} 
Using a level-raising argument (and a result of Larsen on the image of Galois
representations in compatible systems), we prove that for any 
automorphic representation $\pi$ for $\U(3)$, the $l$-adic Galois
representation $\rho_l$ which is attached to $\pi$ by the work of Blasius and 
Rogawski, is the one expected by local Langlands
correspondance at every finite place 
(at least up to semi-simplification and for a density one
set of primes $l$). We rely on the work of Harris and Taylor, who have proved
the same results (for $\U(n)$) assuming the base change of $\pi$ is
square-integrable at one place. As a corollary, every automorphic 
representation which is tempered at an infinite number of places is tempered 
at 
every places.
\end{abstract}

\maketitle
\section{Introduction}

\subsection{Résultats}

Soit $E/F$ une extension CM de corps de nombres, $c$ l'élément non
trivial de $\Gal(E/F)$, $G$ le groupe unitaire à trois variables
défini par $$G(R)=\{g \in \Gl_n(E \otimes_F R),\ {}^tc(g) g=1\}$$ pour toute
$F$-algèbre $R$. On choisit une clôture algébrique ${\bar E}$ de $E$ et 
on pose $\Gamma_E=\Gal(\bar E/E)$.

Dans~\cite{Roglivre}, Rogawski a étudié en grand
détail les représentations automorphes de $G$. Il 
 a construit une application de 
{\it changement de base} attachant à chaque 
représentation automorphe $\pi$ de $G$ une représentation automorphe 
$\pi_E$ de  $G_E = G \times_F E \simeq (\Gl_3)_E$ (voir~\ref{chbaseglobal}). 
De plus, avec Blasius, dans~\cite{br}, il a montré 
l'existence d'un système compatible de  représentations
galoisiennes attaché à $\pi$ (ou à $\pi_E$), que nous explicitons ci-dessous.

Pour chaque place finie $w$ de $E$, choisissons une clôture algébrique 
$\bar E_w$ de $E_w$, et notons $\Gamma_{E_w}:=\Gal(\bar E_w/E_w)$. Choisissons
 également un plongement de $\bar E$ dans $\bar E_w$, ce qui permet 
d'identifier
$\Gamma_{E_w}$ à un sous-groupe {\it de décomposition} $D_w$ 
de $\Gamma_E$ en $w$. 

Pour $l$ un nombre premier, choisissons une clôture algébrique 
$\bar \Q_l$ de $\Q_l$
et notons (à la suite de~\cite[page 6]{ht}) $r_l$ la 
correspondance de
Langlands locale qui associe à une représentation complexe lisse 
irréductible de $\Gl_3(E_w)$ une représentation continue de dimension
$3$ de $\Gamma_{E_w}$ sur $\bar \Q_l$. 

Enfin, pour $r$ une représentation d'un groupe $G$, nous notons
dans cet article $r_{|H}$ la restriction de $r$ à un sous-groupe $H$ de $G$
et $r_{|H}^\ses$ la semi-simplification de cette restriction.

\begin{surtheoreme}[Blasius-Rogawski (voir \cite{br} 
théorème 1.9.1 (a) et théorème 2.2.1.)] \label{thbr}
Soit $\pi$ une représentation automorphe de 
$G$ dont le changement de base $\pi_E$ est
cuspidal. Il existe un corps de nombres $L$, et pour toute place finie
$\mu$ de $L$  une
représentation continue, absolument irréductible,
$\rho_\mu : \Gal(\bar E/E) \rightarrow \Gl_3(L_\mu)$ 
telle que pour toute place finie $v$ de $F$ première à
$N(\mu) \Disc(E/\Q)$ telle que $\pi_v$ est non ramifiée, et pour toute place 
$w$ de $E$ divisant $v$, on ait 
$$(\rho_\mu)_{|D_w}^\ses \simeq r_l((\pi_E)_w).$$

De plus, soit $\pi_E$ est l'induite automorphe d'un caractère de Hecke
$\Psi$ d'une extension cubique $E'$ de $E$, soit les représentations $\rho_\mu$
sont {\it fortement} absolument irréductibles, i.e. la  restriction de
$\rho_\mu$
à  tout sous-groupe ouvert de $\Gal(\bar E/E)$ est absolument irréductible.
\end{surtheoreme}

Si $L$ est un corps de nombres satisfaisant la conclusion du théorème 1,
on dira que $L$ est un {\it corps de définition} de $\pi$. Tout corps
de nombres contenant $L$ est alors aussi un corps de définition de $\pi$.
 
Les $\rho_\mu$ forment donc un système de représentations compatibles
de $\Gal(\bar E/E)$ à coefficients dans $L$. Mais le théorème précédent
ne décrit les
représentations $\rho_\mu$ qu'aux bonnes places $w$ de $E$. 
L'objet de cet article est d'étendre cette description aux autres places.
Nous y parvenons, mais seulement pour un ensemble de densité 1 de places
$\mu$, et à semi-simplification près. Plus précisément, notre résultat principal est le suivant :

\begin{surtheoreme} \label{principal}
Soit $\pi$ une représentation automorphe de $G$ 
dont le changement de base $\pi_E$ est
cuspidal, et $L$ un corps de définition de $\pi$. Il existe un ensemble $S$ 
de densité $1$ de nombres premiers 
tel que pour toute place $\mu$ de $L$ divisant $l \in S$,
et pour toute place finie $w$ de $E$ ne divisant pas $l$,
on ait  $$(\rho_\mu)_{|D_w}^\ses \simeq r_l((\pi_E)_w)^\ses.$$ 
\end{surtheoreme}

Remarquons que si $\pi_E$ est une induite automorphe d'un caractère de
Hecke $\Psi$, le théorème~\ref{principal} découle aisément des
propriétés de compatibilité locale/globale de l'induction automorphe
(cf~\cite{AC}) et du théorème de multiplicité un forte pour $\Gl_3$.
Pour prouver le théorème~\ref{principal}, on peut donc supposer, et l'on
supposera par la suite, que les $\rho_\mu$ sont fortement absolument 
irréductibles.

Le résultat suivant se déduit aisément du théorème~\ref{principal}:
\begin{surcorollaire} 
\begin{itemize}
\item[1)] Si $\pi$ est une représentation automorphe de $G$ dont le changement
de base $\pi_E$ est cuspidal, alors $\pi_E$
est tempérée en toute place $w$ de $E$.
\item[2)] Si $\pi$ est une représentation automorphe de $G$ dont le changement
de base $\pi_E$ est cuspidal, alors $\pi$
est tempérée en toute place $v$ de $F$.
\item[3)] Si $\pi$ est une représentation automorphe auto-duale
de $\Gl_3(F)$ ($F$ un corps totalement réel), tels que $\pi_\infty$ est 
discrète à poids distincts pour toute place réelle $\infty$ de $F$,
alors $\pi$ vérifie la conjecture de Ramanujan.
\end{itemize}
\end{surcorollaire}
\begin{pf} Le premier point se prouve à partir du théorème~\ref{principal}
exactement comme le théorème VIII.1.11 de \cite{ht}. Le second résulte du 
premier. Pour le troisième, pour toute place finie $v$ de $F$, on peut choisir
un corps CM $E/F$ tel que le changement de base $\pi_E$ est cuspidal et 
$v$ décomposé dans $E$. Les hypothèses faites sur $\pi$ assurent alors
que $\pi_E$ provient par changement de base d'une représentation automorphe 
$\pi'$ du groupe unitaire $G$ attaché à $E/F$. Le fait que $\pi_v$ est tempérée
résulte alors de 1).
\end{pf}

\subsection{Méthode}
 
La preuve du théorème~\ref{principal} 
s'appuie sur le résultat essentiel suivant, du à Harris et Taylor :

\begin{surtheoreme}[Harris-Taylor] \label{thht} 
Soit $\pi$ une représentation automorphe de
 $G$ dont le changement de base $\pi_E$ est
cuspidal, et $L$ un corps de définition de $\Pi$. 
{\it On suppose que $(\pi_E)_w$ est de carré intégrable pour au moins une 
place finie $w$ de $E$}. 
Alors pour toute place $\mu$ de $L$ (de caractéristique résiduelle $l$),
et pour toute place finie $w$ de $E$ ne divisant pas $l$,
on a  $$(\rho_\mu)_{|D_w}^\ses \simeq r_l((\pi_E)_w)^\ses.$$ 
\end{surtheoreme}

L'objet du théorème~\ref{principal} est donc 
d'enlever l'hypothèse que $\pi_E$ est de carré intégrable
au théorème~\ref{thht}. Pour cela, on prouve d'abord deux propositions
qui sont des cas particuliers du théorème~\ref{principal} :
\begin{surprop} \label{propinertie} Soit $\pi$ une représentation automorphe de
 $G$ dont le changement de base $\pi_E$ est
cuspidal, et $L$ un corps de définition de $\pi$. 
Il existe un ensemble $S$ de densité $1$ de nombre premiers 
tel que pour toute place $\mu$ de $L$ divisant $l \in S$,
et pour toute place finie $w$ de $E$ ne divisant pas $l$,
on ait  $$(\rho_\mu)_{|I_w}^\ses \simeq r_l((\pi_E)_w)_{|I_w}^\ses,$$ 
où $I_w \subset D_w$ est le sous-groupe d'inertie.
\end{surprop}

On note $G_v$ le groupe $G(F_v)$, et $B_v$ un sous-groupe d'Iwahori de $G_v$.
\begin{surprop} \label{propiwahori} 
Soit $\pi$ une représentation automorphe de $G$ 
dont le changement de base $\pi_E$ est
cuspidal, et $L$ un corps de définition de $\pi$. 
Il existe un ensemble $S$ de densité $1$ de nombre premiers 
tel que pour toute place $\mu$ de $L$ dont la caractéristique résiduelle $l$
appartient à $S$,
et pour toute place finie $v$ de $F$ ne divisant pas $l$, 
si $\pi_v^{B_v} \neq 0$, on ait, pour $w$ place de $E$ au-dessus de $v$, 
  $$(\rho_\mu)_{|D_w}^\ses \simeq r_l((\pi_E)_w)^\ses.$$ 
\end{surprop}

Pour prouver les propositions~\ref{propinertie} et~\ref{propiwahori}, 
l'idée est, en première approximation, de se ramener au 
théorème~\ref{thht}, appliqué à une représentation
automorphe $\pi_n$ qui en vérifie l'hypothèse (de carré intégrable en une
place) et qui
 est congrue à $\pi$ modulo $\mu^n$, pour $n=1,2,\dots$, la 
représentation $\pi_n$ étant obtenue en appliquant le théorème 
d'augmentation du niveau de~\cite{bg} à $\pi$ en une place inerte auxiliaire $v_0$.
Cette méthode consistant à utiliser des augmentations du niveau pour
se ramener au cas où une hypothèse technique est satisfaite a d'ailleurs
déjà été utilisée dans la thèse de Taylor (\cite{taylor}), 
dans un cadre et un but différents
(il s'agissait de montrer l'existence de représentations galoisiennes
attachées aux formes modulaires de Hilbert)

Cependant on bute ici sur deux difficultés : la première est qu'en appliquant
un théorème d'augmentation du niveau à $\pi$ modulo $\mu^n$ avec $n>1$,
 on n'obtient pas,
de représentations automorphes $\pi_n$ mais seulement des formes 
automorphes $f_n$ qu'on ne peut supposer propres pour les opérateurs de Hecke.
On contourne ce problème en imposant à $f_n$ d'apparaître dans un
espace de formes automorphes admettant un type de Bushnell-Kutzko 
adéquatement choisi (pour prouver la proposition 1) ou d'être propre modulo
$\mu^n$ pour un caractère adéquat du centre de l'algèbre de Hecke-Iwahori 
(pour prouver la proposition 2). 

La seconde difficulté consiste à montrer
l'existence de bonnes places inertes auxiliaires : 
on ne peut pour cela utiliser le  théorème 2 de~\cite{bg}, qui montre 
l'existence de places auxiliaires $v_0$, mais telles que  $\mu$ est non 
{\it normale} pour $v_0$ (voir~\ref{banal}), ce qui empêche de montrer que la forme obtenue 
$f_n$ vérifie les hypothèses du théorème~\ref{thht}. On utilise à la place 
un résultat
récent de Larsen sur l'image des représentations galoisiennes dans un système
compatible (\cite{larsen}), qui combiné avec un résultat de Steinberg,
permet de montrer l'existence de bonnes places inertes $v_0$, telles que
$\mu$ est normal pour $v_0$ (quoique non banal),
où augmenter le niveau, à condition d'exclure un ensemble de mesure nulle de nombres premiers $l$.

Enfin, pour prouver le théorème~\ref{principal}, on utilise 
un changement de base non galoisien pour ramener le calcul
de la trace d'un élément de $D_w-I_w$ dans le cas général 
à la proposition 2, ce qui, combiné avec la proposition $1$ 
(déterminant la trace des actions des éléments de $I_w$) permet de conclure. 

\section{Notations et rappels} 

\subsection{Notations}
On garde les notations de l'introduction, $E/F$ est une extension $CM$,
$c$ l'élément non trivial de $\Gal(E/F)$. On choisit un relevé $\gamma$ de $c$
dans $\Gamma_F = \Gal(\bar F/F)$, tel que $\gamma^2=1$. 

Soit $G$ le groupe unitaire défini dans l'introduction. On notera $G_v$ le
groupe $G(F_v)$. On notera $q_v$ le cardinal résiduel de $F_v$, et 
$\pig_v$ une
uniformisante de $F_v$. 

\subsection{Le changement de base de Rogawski}

\subsubsection{Changement de base global}
\label{chbaseglobal}

Rogawski définit une partition des représenta\-tions automorphes de $G$
en $A$-paquets globaux, et pour chaque $A$-paquet $\Pi$, son changement de base
$\pi_E$ qui est une représentation automorphe de $G_E$. On peut donc 
définir sans ambiguïté le {\it changement de base d'une représentation 
automorphe} $\pi$ de $E$ comme le changement de base du $A$-paquet auquel elle 
appartient.

\subsubsection{Changement de base local}
\label{chbaselocal}

Soit $v$ une place de $F$.
Rogawski définit, suivant les conjectures d'Arthur, des parties finies de 
l'ensemble des classes d'isomorphismes
de représentations lisses complexes irréductibles de $G_v$ qu'il appelle
$A$-paquet locaux. Si $v$ est décomposée, tous les $A$-paquets sont
des singletons.

Soit $w$ une place de $E$ au-dessus de $F$. Si $\Pi_v$ est un $A$-paquet
local, on peut définir son changement de base (local) à $E_w$ 
qui est une représentation lisse complexe irréductible de $\Gl_3(E_w)$.  
Les $A$-paquets globaux sont produits tensoriels de $A$-paquets locaux
et le changement de base global est compatible au changement de base local.

Une représentation complexe lisse irréductible de $G_v$ peut appartenir
à plusieurs $A$-paquets locaux. Cependant :
\begin{lemme} \label{lemmechbaselocal} 
Si deux $A$-paquets $\Pi_v$ et $\Pi'_v$ contiennent une représentation en commun, leurs changements de bases sont des sous-quotients d'une même
induite parabolique. En particulier, les images de ces deux changements de 
base par la correspondance de Langlands ont même semi-simplification.

De plus, si $\Pi_v$ et $\Pi'_v$ contiennent deux représentations $\pi_v$ et $\pi'_v$
ayant même support supercuspidal à équivalence inertielle près, il en va de 
même de leur changement de base. 
\end{lemme}
\begin{pf}
Il n'y a rien à prouver pour $v$ décomposée. On suppose $v$ inerte ou ramifiée.

D'après \cite{Roglivre},
les seules représentations appartenant à plusieurs $A$-paquets
sont celles notées  $\pi^s(\chi)$ {\it loc. cit., page 199}, 
pour $\chi$ un caractère de $\U(1)$. Ces 
représentations appartiennent exactement à deux $A$-paquets, 
$\{\pi^s(\xi),\pi^2(\xi)\}$
(qui est tempérée), et $\{\pi^s(\xi),\pi^n(\xi)\}$ (qui ne l'est pas).
Les changements de base de ces deux $A$-paquets sont respectivement la
sous-représentation irréductible et la représentation quotient de l'induite
du Borel du caractère $(\xi_E,\xi_E|\ |,1)$. 

Le ``en particulier'' résulte par exemple de \cite{he}.

Le ``de plus'' résulte de \cite[prop. 13.2.2]{Roglivre}.
\end{pf}

Si $\pi_v$ est la composante locale d'une représentation automorphe, elle appartient à au moins un $A$-paquet local. 
La notation $r_l((\pi_v)_{E_w})^\ses$ est donc définie sans ambiguïté par
$r_l((\Pi_v)_{E_w})^\ses$, où $\Pi_v$ est n'importe quel paquet local contenant
$\pi_v$.

\subsection{Rappels locaux}

\subsubsection{Algèbres de Hecke non ramifiées}

Si $v$ est inerte ou ramifiée, $G_v$ est un groupe unitaire quasi-déployé 
de rang un. Si $v$ est inerte, $G_v$ est de plus non ramifiée. On dispose
alors d'une classe de conjugaison privilégiée de sous-groupes compacts 
maximaux de $G_v$, la classe de conjugaison hyperspéciale. Si $K_v$ appartient
à cette place l'algèbre de Hecke $\Hecke(G_v,K_v)$ des fonctions à support 
compact, $K_v$-invariantes à gauche et à droite, à valeurs dans $\Z$, munies
de la convolution, est commutative, isomorphe à $\Z[T_v]$. Ici, $T_v$ est 
l'opérateur de Hecke standard, défini comme la fonction caractéristique de la
double classe de la matrice diagonale $(\pig_v,1,\pig_v^{-1})$. 
On le notera aussi parfois  $T_w$, où $w$ 
est la place de $E$ au-dessus de $v$

Si $v$ est une place décomposée, le choix d'une place $w$ de $E$ définit
un isomorphisme $G_v \simeq \Gl_3(F_v) \simeq \Gl_3(E_w)$, canonique à 
automorphismes intérieurs près. Si $K_v$ est le compact maximal correspondant à $\Gl_3(\anneau_v)$ via cet isomorphisme, on notera $T_w$ l'opérateur de Hecke
dans $\Hecke(G_v,K_v)$ donnée par la fonction caractéristique de la double 
classe de la matrice diagonale $(\pig_v,1,1)$. 

\subsubsection{}

\label{br2}

On peut reformuler le théorème~\ref{thbr} en termes des opérateurs de Hecke
ainsi définie. En gardant les notations de ce théorème, pour toute place $v$ 
telle que $\pi_v$ est non ramifiée, $w$ place de $E$ au-dessus de $v$,
si $\lambda_w$ est  la valeur propre par laquelle agit l'opérateur $T_w$ sur 
la droite  $\pi_v^{K_v}$, alors 
$$\lambda_w = \tr \rho_\mu(\Frob_w).$$

\subsubsection{Caractéristiques normales}
\label{banal}

Si $v$ est inerte, on dit que $l$ est une caractéristique {\it normale} 
(resp. {\it banales}) pour $G_v$ (ou pour $v$) si $l$ ne divise pas 
$q_v(q_v^3+1)$ (resp. $q_v(q_v-1)(q_v^3+1)$).

\subsubsection{Centre de Bernstein}

Nous rappelons un fragment de la théorie du centre de Bernstein.
Soit $B_v$ un sous-groupe d'Iwahori de $G_v$, $K_v$ un compact maximal. On 
note $\ZZ(G_v,B_v)$ le centre de l'algèbre de Hecke-Iwahori 
$\Hecke(G_v,B_v)$. Alors
\begin{lemme} \label{bernstein} Il existe un isomorphisme d'algèbre
 $$b :\ \ZZ(G_v,B_v)\otimes \C \simeq  \Hecke(G_v,K_v) \otimes \C,$$ tel que, 
pour toute 
représentation $\pi$ lisse irréductible de $G_v$ telle que 
$\pi^{B_v} \neq 0$,
\begin{itemize} \item[1)] $\ZZ(G_v,B_v)$ agit sur $\pi^{B_v}$ par un caractère 
$\psi$,
et $\Hecke(G_v,K_v)$ agit sur $\pi^{K_v}$ (éventuel\-lement nul) par le
caractère $\psi \circ b^{-1}$.
\item[2)] Il existe une unique représentation $\pi'$ avec 
${\pi'}^{K_v} \neq 0$, et telle que $\Hecke(G_v,K_v)$ agisse sur 
$\pi^{K_v}$ par $\psi \circ b^{-1}$. De plus $\pi$ et $\pi'$ apparaissent 
dans une même induite non ramifiée.
\item[3)] La représentation galoisienne 
$r_l((\pi)_{E_w})^\ses$ est non ramifiée, et $\tr r_l((\pi)_{E_w})(\Frob_w) =
\psi \circ b^{-1} (T_w).$
\end{itemize}
\end{lemme} 
\begin{pf} D'après la théorie classique des représentations non ramifiées,
toute représentation lisse irréductible $\pi$ apparaît comme facteur de
Jordan-Hölder d'une induite parabolique $I$ d'un caractère non ramifié, 
qui contient un unique facteur non ramifiée $\pi'$, et qu'on peut supposer
indécomposable. D'après l'équivalence de catégorie bien connue de Borel,
$I^{B_v}$ est un module indécomposable sur $\Hecke(G_v,B_v)$. Le centre 
$\ZZ(G_v,B_v)$ agit donc par un caractère $\psi$ sur $I^{B_v}$, donc sur
$\pi^{B_v}$ et $\pi'^{B_v}$. L'existence de l'isomorphisme $b$ vérifiant 1)
résulte de~\cite{bernstein}; 2) résulte alors de ce qui précède, et 3) de
\cite{he}.
\end{pf} 
\section{Existence de places inertes et normales où augmenter le niveau}
\subsection{Énoncé}
\subsubsection{}
\label{explaces}
Dans toute cette partie, $\pi$ est une 
représentation automorphe pour $G$ dont le changement de base à $E$,
est cuspidal, $w$ est une place de $E$ et $v$ est la place de $F$
divisant $E$. On choisit un corps de définition $L$
de $\pi$, et pour $\mu$ place finie de $L$, on note $\rho_\mu$ la 
représentation galoisienne associée, qu'on suppose fortement absolument
irréductible.

Le but de cette sous-partie est de prouver le résultat suivant : 
\begin{prop} \label{exbonnesplaces} 
Il existe un ensemble $S$ de nombre premiers, de
densité $1$, tel que pour toute place finie $\mu$ de $L$ divisant
 une place de $S$, et pour tout entier $n \geq 1$, il existe
une infinité de places $v_0$ de $F$, inertes dans $E$, telles que $\pi_{v_0}$ 
est
non ramifiée, et telles que la valeur propre $\lambda=\lambda_{v_0}$ 
de l'opérateur de Hecke  standard $T_{v_0}$ sur $\pi_{v_0}$, et le cardinal résiduel $q=q_{v_0}$ de $F_{v_0}$ vérifient 
\begin{eqnarray} \lambda &\equiv & q(q^3+1)
 \pmod{\mu^n} \\
q &\equiv&  1 \pmod{\mu^n}.
\end{eqnarray}
\end{prop}

\subsection{Réduction à une propriété de l'image de $\rho_\mu$}

\subsubsection{}
D'après~\cite[6.2.8]{bg}, il existe une base de $L_\mu^3$ telle que dans
cette base, la représentation $\rho_\mu$ définisse un morphisme
$\rho_\mu: \Gal(\bar E/E) \rightarrow \Gl_3(\anneau_\mu)$,
vérifiant pour tout $g \in G_E$
\begin{eqnarray} \label{antiauto} \rho_\mu(\gamma g \gamma^{-1}) = A {}^t
\rho_\mu(g)^{-1} A^{-1},
\end{eqnarray}
où $A \in \Gl_3(\anneau_\mu)$. Comme $\gamma$ est une involution, 
on voit en appliquant deux fois~(\ref{antiauto}) que $A {}^t A^{-1} $ est dans 
le commutant de 
$\rho_\mu$, donc est un scalaire $\alpha$. D'où ${}^t A = \alpha A$, 
ce qui implique $\alpha = \pm 1$, et le cas antisymétrique $\alpha=-1$ 
est impossible car $A$ est inversible. On a donc 
\begin{eqnarray} 
\label{asymetrique} {}^t A=A. \end{eqnarray}

\subsubsection{}
Notons $\rho_n$ la réduction de $\rho_\mu$ modulo $\mu^n$, et $\omega_{l^n} : \Gal(\bar F/F) \rightarrow (\Z/l^n\Z)^\ast$ le
caractère cyclotomique. La proposition~\ref{exbonnesplaces} se réduit à la proposition
suivante : 
\begin{prop} \label{image} Il existe $g \in \Gal(\bar E/E)$ tel que
$\rho_n(g) = a \Id$ avec $a \in (\anneau/\mu^n)^\ast$ et $\omega_{l^n}(g)=-1$.
\end{prop}

\subsubsection{}
Prouvons que le proposition~\ref{image} implique la 
proposition~\ref{exbonnesplaces}. 
$\tilde G$ le produit semi-direct $$\Gl_3(\anneau/\mu^n) \sd C,$$ où
$C$ est le groupe à deux éléments $\{1,c\}$ et $c$ agit sur $G$ par $g
\mapsto {}^t g^{-1}$. Les relations~(\ref{antiauto},
\ref{asymetrique})) permettent  de prolonger $\rho_n$ en un morphisme, 
noté $\tilde \rho_n : \Gal(\bar E/F) \rightarrow \tilde G$ par
\begin{eqnarray}
\tilde \rho_n(\sigma) &=& \rho(\sigma) \sd 1, \ \sigma \in \Gal(\bar E/F) \\
\tilde \rho_n(\gamma) & = & A \sd c 
\end{eqnarray}

Si $g$ est comme dans l'énoncé de la proposition, on a $\tilde \rho_n
(g \gamma)
= (aA) \sd c$ et $\omega_{l^n}(g \gamma)=1$. Il existe donc une infinité
de places $v_0$ de $F$ telles que $\tilde \rho_n
(\Frob_{v_0}) = (a A) \sd c$ et $\omega_{l^n}(\Frob_{v_0})=1$. Ces places ne
peuvent être décomposées, un infinité d'entre elles sont donc
inertes. Pour $v_0$ une telle place, on a $q=1$, et si $w_0$ est la place
de $E$ au-dessus de $v_0$, 
$$\rho_n(\Frob_{w_0})=\tilde \rho_n(\Frob_{v_0})^2= (a A) {}^t({a A})^{-1}=\Id$$
(cf. \cite[6.2.8]{bg}),
ce qui implique que $\lambda \equiv q(q^3+1) \pmod{\mu^n}$ d'après la
transformation de Satake, cf. \cite[3.7.1]{bg}

\subsubsection{} Le reste ce cette partie est consacrée à la preuve
de la proposition~\ref{image}. Nous distinguerons deux cas (i) et (ii), suivant
la disjonction du théorème 2.2.1 de \cite{br}, point (b) : le cas (i) 
est celui où les représentations $\rho_\mu$ sont fortement absolument
irréductibles, i.e. sont absolument irréductibles  et le restent après 
restriction  à tout sous-groupe d'indice
fini du groupe $\Gal(\bar E/E)$. Le cas (ii) est celui des représentations qui sont des induites
d'un sous-groupe d'indice $3$.

\subsection{Preuve de la proposition~$\ref{image}$}

\subsubsection{}
Pour tout nombre premier $l$, et toute place $\mu$ de $L$ au-dessus de
$l$, $\rho_\mu$ est une représentation de $\Gal(\bar E/E)$ de
dimension 3 sur $L_\mu$ qu'on peut 
également voir comme une représentation de dimension $3[L_\mu:\Q_l]$
sur $\Q_l$. Posons 
$$\rho_l=\oplus_{\mu | l} \rho_\mu$$ qui est une
représentation de dimension $$3\sum_{\mu | l}[L_\mu:\Q_l]=3[L:\Q].$$
Les $\rho_l$ forment un système compatible de représentations
$l$-adiques.

Notons $\Gamma_l$ l'image de $\rho_l$  (dans $\Gl_{3[L:\Q]}(\Q_l)$)
et $G_l$ son adhérence Zariski. Comme les représentations $\rho_\mu$ (sur
$L_\mu$) sont irréductibles, elles sont semi-simples en tant que
représentation sur $\Q_l$, ainsi que les $\rho_l$. Les groupes $G_l$
sont donc réductifs. On note $G_l^0$ la
composante neutre de $G_l$, auquel on peut appliquer un résultat 
de Serre (cf. \cite[théorème, page 15]{Serre}): 

\begin{lemme} Il existe une extension finie $E'$ de $E$, tel que pour tout $l$,
$$\rho_l(\Gal(\bar E/E'))=\Gamma_l \cap G_l^0(\Q_l).$$ \end{lemme}

Quitte à remplacer $E$ par $E'$, on peut donc supposer que $G_l=G_l^0$.
Comme $\rho_\mu$ est fortement absolument irréductible, 
Ce changement de $E$ en $E'$ n'affecte pas 
l'hypothèse que $\rho_\mu$ est absolument irréductible.

\subsubsection{}
On note, suivant~\cite{larsen}, $G_l^\ad$ le groupe adjoint de $G_l$ et
$G_l^\ssc$ le revêtement universel de $G_l^\ad$. Ce sont des groupes
réductifs connexes, et $G_l^\ssc$ et $G_l^\ad$ sont même semi-simples.

D'après~\cite{lp}, il existe un ensemble de densité 1 de
nombre  premiers, tel que pour tout nombre premier $l$ de cet ensemble, 
le groupe réductif $G_l^\ssc$ est non ramifié (cf. \cite{tits}). 
On peut donc parler de sous-groupes compacts maximaux hyperspéciaux 
de $G_l^\ssc(\Q_l)$ pour $l \not \in S$ 
(cf. \cite{tits}). 
Notons enfin $\Gamma_l^\ad$ l'image de $\Gamma_l$ dans
$G_l^\ad(\Q_l)$ par l'application canonique, et $\Gamma_l^\ssc$ l'image
réciproque de $\Gamma_l^\ad$ dans $G_l^\ssc(\Q_l)$.

$$\xymatrix{\Gamma_l \ar[rd] \ar[rr]  & & \Gamma_l^\ad \ar[rd] & & \Gamma_l^\ssc \ar[ll]
\ar[rd] & \\
& G_l(\Q_l)  \ar[rr] & & G_l^\ad(\Q_l) & & G_l^\ssc(\Q_l) \ar[ll]}$$ 

Le théorème suivant est le résultat principal de~\cite{larsen}.
\begin{theoreme}[Larsen] Il existe un ensemble de densité $1$ 
de nombres premiers tel que pour tout $l$  de cet ensemble,
$\Gamma_l^\ssc$ est un sous-groupe compact maximal hyperspécial 
de $G_l^\ssc(\Q_l)$.
\end{theoreme}

\begin{corollaire} \label{abtriv} Il existe un ensemble $S$ de densité $1$
de nombres premiers tel que pour tout $l$ de cet ensemble,
$\Gamma_l^\ssc$ n'a pas de quotient d'ordre $2$.
\end{corollaire}
\begin{pf} On prend pour ensemble $S$ de nombre premiers celui du théorème 
précédent privé de  $2$. Pour $l$ en
dehors de cet ensemble, $G_l^\ssc$ a un (unique) modèle à fibres semi-simples
connexes sur $\Z_l$ tel que $\Gamma_l^\ssc=G_l^\ssc(\Z_l)$. Comme le
noyau de la réduction, surjective d'après le lemme de Hensel, 
$G_l^\ssc(\Z_l) \rightarrow G_l^\ssc(\F_l)$ est un pro-$l$-groupe avec $l
\neq 2$, il suffit de voir que $G_l^\ssc(\F_l)$ n'a pas de quotient d'ordre
2. Or d'après un théorème de Steinberg (\cite{steinberg}, voir 
aussi~\cite[page 406]{pr}), comme  le $\F_l$-groupe algébrique $G_l^\ssc$
est  semi-simple connexe et simplement connexe,
 le groupe de ses points rationnels $G_l^\ssc(\F_l)$ est d'abélianisé trivial.
\end{pf}

Dorénavant, $S$ désigne un ensemble de nombres premiers 
comme dans le corollaire ci-dessus.

\subsubsection{} 
Notons maintenant $\Gamma_\mu$ l'image de $\rho_\mu$ dans
$\Gl_3(L_\mu)$, et $G_\mu$ l'adhérence Zariski de $\Gamma_\mu$ dans 
$\Gl_3(L_\mu)$ considéré comme groupe algébrique sur $\Q_l$,
i.e. dans $\Res^{L_\mu}_{\Q_l} \Gl_3$. 

La projection $p:\  \rho_l \rightarrow \rho_\mu$ définit un morphisme
 surjectif $p$ de $\Gamma_l$ sur $\Gamma_\mu$ qui se prolonge en un morphisme
surjectif (noté aussi $p$) de $\Q_l$-groupes algébriques $G_l \rightarrow G_\mu$. Le groupe $G_\mu$ est donc connexe. On définit $G_\mu^\ad$ comme 
le  groupe adjoint de $G_\mu$, et $G_\mu^\ssc$ comme le revêtement 
universel  de $G_\mu^\ad$. On définit  $\Gamma_\mu^\ad$ comme l'image de 
$\Gamma_\mu$,  et $\Gamma_\mu^\ssc$ comme l'image réciproque de 
$\Gamma_\mu^\ad$.

On a donc le diagramme commutatif suivant :
$$\xymatrix{\Gamma_l \ar[rd] \ar[rr] \ar[dd] & & \Gamma_l^\ad \ar[rd] \ar[dd]& & \Gamma_l^\ssc \ar[ll]
\ar[rd]\ar[dd] & \\ 
& G_l(\Q_l)  \ar[rr]\ar[dd] & & G_l^\ad(\Q_l) \ar[dd] & & G_l^\ssc(\Q_l) \ar[ll] \ar[dd]
\\ \Gamma_\mu \ar[rd] \ar[rr]  & & \Gamma_\mu^\ad \ar[rd] 
& & \Gamma_\mu^\ssc \ar[ll]^{\ \ \ \ \ \tau}
\ar[rd] & \\
& G_\mu(\Q_l)  \ar[rr] & & G_\mu^\ad(\Q_l) & & G_\mu^\ssc(\Q_l) \ar[ll]^{\tau}}
$$ 
où les flèches verticales se déduisent de $p$ par fonctorialité et sont toutes
surjectives par construction.

\begin{lemme} \label{gammasc}
Si $\mu | l$ et $l \in S$, $\Gamma_\mu^\ssc$ n'a pas de 
quotient d'ordre $2$.
\end{lemme}
\begin{pf} Comme c'est un quotient de $\Gamma_l^\ssc$, cela résulte du
corollaire~\ref{abtriv}. \end{pf}
\begin{lemme}
Le groupe $G_\mu^\ad(\Q_l)/\tau(G_\mu^\ssc(\Q_l))$ est abélien de $3$-torsion.
\end{lemme} 
\begin{pf} Considérons $G_\mu$ comme un sous-groupe algébrique de 
$$\Res_{L_\mu/\Q_l} \Gl_3.$$ 
Le groupe $G_\mu \times \bar \Q_l$ 
est donc un sous-groupe
e algébrique de  
$$(\Res_{L_\mu/\Q_l} \Gl_3) \times \bar \Q_l=\prod_W \Gl_3$$ 
où $W$ est l'ensemble des plongements de $L_\mu$ dans $\bar \Q_l$. 
On en déduit que le
groupe $G^\ad \times \bar \Q_l$ est un sous-groupe de $\prod_W \Sl_3$.
Soit $Z$ son centre : comme $\rho_\mu$ est absolument irréductible,
c'est un sous-groupe du centre de $\prod_W \Sl_3$, et il est donc
de $3$-torsion.

Considérons la suite exacte courte de groupes munis d'une action 
de $\Gal(\bar \Q_l/L_\mu)$ : 
$$0 \rightarrow Z(\bar \Q_l) \rightarrow G_\mu^\ssc(\bar \Q_l) 
\stackrel{\tau}\rightarrow
G_\mu^\ad(\bar\Q_l) \rightarrow 0.$$
La suite exacte longue de cohomologie associée identifie 
$G_\mu^\ad(\Q_l)/\tau(G_\mu^\ssc(\Q_l))$ à un sous-groupe de 
$H^1(L_\mu, Z(\bar \Q_l))
$ qui est abélien de 
$3$-torsion. 
\end{pf}
\begin{lemme} Si $\mu | l$ et $l \in S$, $\Gamma_\mu^\ad$ n'a pas de 
quotient d'ordre $2$. 
\end{lemme}
\begin{pf} Le groupe $G_\mu^\ad(\Q_l)/\tau(G_\mu^\ssc(\Q_l))$ 
est abélien de $3$-torsion, d'après le lemme précédent. Il en va donc de même 
du groupe 
$\Gamma_\mu^\ad /\tau(\Gamma_\mu^\ssc)$, qui n'a donc pas de quotient 
d'ordre $2$. 
D'après le lemme~\ref{gammasc}, $\tau(\Gamma_\mu^\ssc)$ n'a pas non 
plus de quotient d'ordre $2$. Il en va donc de même de $\Gamma_\mu^\ad$.
\end{pf}
\begin{lemme} Le centre $Z(\Gamma_\mu)$  est composé d'homothéties de 
$L_\mu^\ast$, et $\Gamma_\mu/Z(\Gamma_\mu)$ n'a pas de quotient d'ordre $2$.
\end{lemme}
\begin{pf} La première assertion résulte du lemme de Schur, qui montre aussi
que $Z(\Gamma_\mu) = Z(G_\mu(\Q_l)) \cap \Gamma_\mu$  et la seconde résulte 
alors du lemme précédent, car 
$\Gamma_\mu/Z(\Gamma_\mu)=\Gamma_\mu/(Z(G_\mu(\Q_l))\cap \Gamma_\mu) = 
\Gamma_\mu^\ad$.
\end{pf}

\subsubsection{Fin de la preuve de la proposition~$\ref{image}$}

Considérons le caractère $\omega_{l^n}^{(l^n-1)/2)}$. 
Son image est $\{\pm 1\}$. D'après le lemme de Goursat et le lemme précédent,
l'application 
$\rho_\mu \times \omega_{l^n}^{(l^n-1)/2)}:\ Gal(\bar E/E) \rightarrow 
\Gamma_\mu/Z(\Gamma_\mu)$ est surjective. Il y a donc un 
$g'\in \Gal(\bar E/E)$ dont l'image est $(1,-1)$. Posons $g=g'^{(l^n-1)/2}$.
Alors $\rho_\mu(g) \in Z(\Gamma_\mu)$ donc est un scalaire, et 
il en va de même de sa réduction $\rho_n(g)$. De plus $\omega_{l^n}(g)=-1$,
ce qui prouve la proposition.

\section{Preuve de la proposition~$\ref{propinertie}$}

On reprend les mêmes notations que dans la partie~\ref{explaces}.

\subsection{Choix d'un niveau et d'un type pour la preuve}

\begin{lemme} \label{type}
Il existe un sous-groupe compact ouvert $K_v$ de $G_v$ et une 
représentation $J_v$ de $K_v$ tel que 
\begin{itemize}
\item[a.] $\Hom_{K_v}(J_v,\pi_v) \neq 0$
\item[b.] Pour $\tau$ une représentation lisse irréductible de $G_v$, si 
$\Hom_{K_v}(J_v,\pi_v)$ on a pour toute place $\mu$ de $L$ 
$$r_l(\tau_{E_w})^{\sss}_{|I_w} \simeq r_l(\pi_{E_w})^\sss_{|I_w}.$$
\end{itemize}
\end{lemme}
\begin{pf}
D'après~\cite[page 772]{bus} si $v$ est décomposé, 
ou le théorème principal de~\cite{blas} si $v$ est inerte ou 
ramifiée, il existe un type $(K_v,J_v)$, tel que pour toute représentation 
lisse irréductible $\tau$ de $G_v$, $\Hom_{K_v}(J_v,\tau) \neq 0$ si et seulement si
$\tau$ et $\pi_v$ ont même support supercuspidal à équivalence inertielle
près. Pour un tel type, la propriété a. est évidemment satisfaite.

Si $v$ est décomposé, on a $G_v=\Gl_3(E_w)$ ; soit 
$\{ \pi_1, \dots, \pi_k\}$ le
support supercuspidal de $\pi$ (on a $1 \leq k \leq 3$, et les $\pi_i$ sont des représentations supercuspidales irréductibles de $\Gl_{n_i}(E_w)$, avec
$\sum_{i=1}^k {n_i} =3$).
Si $\Hom_{K_v}(J_v,\tau) \neq 0$, le support supercuspidal de $\tau$ est
donc $\{ \pi_1 \otimes \chi_1 \circ \det, \dots, \pi_k\otimes 
\chi_k\circ\det\}$,
où les $\chi_i$ sont des caractères non ramifiés de $E_w^\ast$.
Les propriétés de compatibilité à 
l'induction parabolique (\cite[propriété 2, page 6]{ht}) et au twist par un 
caractère (\cite[propriété 3, page 6]{ht}) de $r_l$, et le fait que deux constituants d'une même induite parabolique aient, à semi-simplification près, la même image par $r_l$ (cf~\cite{he}) implique que 
\begin{eqnarray} r_l(\pi_v)^\ses & \simeq & \oplus_{i=1}^k r_l(\pi_i)\\
r_l(\tau)^\ses  & \simeq & \oplus_{i=1}^k r_l(\pi_i) \otimes r_l(\chi_i)
\end{eqnarray}
ce qui implique 
$(r_l(\tau)^\ses)_{|I_w} \simeq (r_l(\pi_v)^\ses)_{|I_w}$.
Comme la restriction à $I_w$ commute à la semi-simplification pour une 
représentation de $D_w$ (cela résulte de~\cite[4.2.1]{Tate}), 
on a donc
$$r_l(\tau)^\ses_{|I_w} \simeq r_l(\pi_v)^\ses_{|I_w}.$$

Supposons $v$ inerte ou ramifiée dans $E$. Soit $\tau$ tel que $\Hom_{K_v}(J_v,\tau) \neq 0$. Soit $\tau_{E_w}$ et $(\pi_v)_{E_w}$ les changements de base
de n'importe quel $A$-paquet contenant $\tau$ et $\pi_v$. 
D'après le lemme~\ref{lemmechbaselocal}, $\tau_{E_w}$ et $(\pi_v)_{E_w}$ 
ont même support supercuspidal à équivalence inertielle près. On est alors ramené au cas précédent.
\end{pf}

\subsubsection{} Soit $K_v,J_v$ comme dans le lemme. 
Soit $K^v=\prod_{v'\neq v} K_{v'}$
 un sous-groupe compact ouvert de $\prod_{v'\neq v} G(F_{v'})$ tel que 
\begin{eqnarray} \label{kvinv} \pi^{K^v} \neq 0.\end{eqnarray}

Posons $K=K_vK^v$ et $J=J_v\otimes 1$.

 Pour toute $\anneau_\mu$-algèbre $R$,
on note $S_{K,J,\pi_\infty,R}$ le $R$-module (libre de rang fini) 
des formes automorphes
de niveau $K$, type $J$ et même poids que $\pi$ qui sont définis sur $R$,
défini en~\cite[5.2]{bg}.

\subsection{Formes anciennes et nouvelles}

\subsubsection{} \label{v0}
Soit $v_0$ une place inerte de $F$, distincte de $v$, et telle que 
$K_{v_0}$ est hyperspécial. On note $T=T_{v_0}$ l'opérateur de Hecke standard
de $\Hecke(G_{v_0},K_{v_0})$ et $q$ le cardinal résiduel de $F_{v_0}$.
Soit $B_{v_0}$ un sous-groupe d'Iwahori de $G_{v_0}$  contenant $K_{v_0}$
et soit $B:=B_{v_0} K^{v_0}$ où $K^{v_0}:=\prod_{v'\neq v_0} K_{v'}$. 

On note $N_{B,J,\pi_\infty,R}$ et $O_{B,J,\pi_\infty,R}$
les $R$-modules (libres de rang finis)
 des formes nouvelles et anciennes (respectivement) en $v_0$ sur $R$ 
de niveau $B$, de type $J$ et de poids $\pi_\infty$,
tels qu'ils sont définis en~\cite[5.3.6]{bg}.
La formation de ces espaces commute à tout changement de base 
$R \rightarrow R'$.

\subsubsection{} \label{actionhecke}
On rappelle que si $\Sigma$ est l'ensemble fini des places 
$v'$ telles que $K_{v'}$
n'est pas hyperspécial, l'algèbre de Hecke commutative 
$$\Hecke^\Sigma=\otimes_{v'\not \in \Sigma}\Hecke(G_{v'},K_{v'})$$ agit par 
opérateurs $R$-linéaire sur les espaces $S_{K,J,\pi_\infty,R}$, 
$N_{B,J,\pi_\infty,R}$ et  $O_{B,J,\pi_\infty,R}$, et ceci fonctoriellement 
en $R$. 
(cf.~\cite[5.2.2 et 5.3.4]{bg}). Rappelons en particulier
que l'action de $\Hecke(G_{v_0},K_{v_0})=\Z[T]$ sur les espaces $N_{B,J,\pi_\infty,R}$ et  $O_{B,J,\pi_\infty,R}$ est définie en identifiant 
$\Hecke(G_{v_0},K_{v_0})$ 
à une sous-algèbre de $\Hecke(G_{v_0},B_{v_0})$,
l'opérateur $T$ étant identifié à l'opérateur $T_B$ défini en \cite[3.3.2]{bg}.
Aux autres places $v'\not \in \Sigma$, l'action de
 $\Hecke(G_{v'},K_{v'})$ est l'action évidente.

Rappelons aussi (\cite[5.1.3]{bg}) 
que si $R$ est une $L_\mu$-algèbre, les actions de 
$\Hecke^\Sigma$ sur les espaces considérés sont semi-simples


\subsubsection{} \label{N+}
Pour tout $\anneau_\mu$-algèbre $R$,
notons $N^+_{B,J,\pi_\infty,R}$ (resp. $N^-_{B,J,\rho,R}$) 
le noyau $\Ker(T-q(q^3+1))$ (resp. $\Ker (T + (q^3+1))$) 
sur $N_{B,J,\rho,R}$. 
\begin{lemme} \label{Nplus} Supposons que $(q^3+1) \not \equiv 0 \pmod{\mu}$.
On a $$N_{B,J,\pi_\infty,R}=N_{B,J,\pi_\infty,R}^+ \oplus 
N_{B,J,\pi_\infty,R}^{-}$$ et pour tout morphisme 
$R \rightarrow R'$ de $\anneau_\mu$-algèbres on a
$$N_{B,J,\pi_\infty,R'}^{\bullet} = N_{B,J,\pi_\infty,R}^{\bullet} 
\otimes_R R'$$ 
pour $\bullet \in \{+,-\}$ .
\end{lemme}
\begin{pf} Comme $q(q^3+1) \neq -(q^3+1) \pmod{\mu}$ par hypothèse, il suffit 
de montrer que l'opérateur $T$ (ou $T_B$)  sur $N^{B,J,\pi_\infty,R}$ est 
annulé par le polynôme $(X-q(q^3+1))(X+(q^3+1))$. Par fonctorialité il suffit 
de le faire pour $R=\anneau_\mu$, et donc comme 
$N_{B,J,\pi_\infty,\anneau_\mu} \subset N_{B,J,\pi_\infty,\bar \Q_l}$
pour $R=\bar \Q_l \simeq \C$. Comme $T$ agit de manière 
semi-simple sur cet espace, il suffit de voir que pour 
$f \in N_{B,J,\pi_\infty,\C}$ propre pour $\Hecke^\Sigma$, la valeur propre 
$\lambda$ de
$T$ est $-(q^3+1)$ ou $q(q^3+1)$. 
Or d'après\cite[6.1.2]{bg}, à une telle forme nouvelle propre $f$ est 
attachée une représentation automorphe $\pi$, telle que $\pi_{v_0}$ est soit la Steinberg $\st$, soit la représentation $\pi^s$ (voir~\cite[3.6.5]{bg}) 
et que $T$ agisse par $\lambda$ sur $\pi_{v_0}^{B_{v_0}}$.
 Le lemme résulte donc de~\cite[6.7.2]{bg}
\end{pf}

\subsection{Augmentation du niveau}

\subsubsection{}
D'après le point a. du lemme~\ref{type} et~(\ref{kvinv}), la 
représentation automorphe $\pi$ définit un élément de $S_{K,J,\pi_\infty,\anneau_\mu}$,
propre pour $\Hecke^\Sigma$ de caractère propre noté $\psi$ 
(cf. \cite[6.7.1]{bg}) 

\begin{prop} Pour tout entier $n \geq 1$, on peut choisir une place $v_0$ 
de $F$ comme en~\ref{v0}, et tel qu'il existe
un vecteur non divisible $f \in N^+_{B,J,\pi_\infty,\anneau_\mu}$, vérifiant
$$U f \equiv \psi(U) f \pmod{\mu^n},$$
pour tout $U \in \Hecke^\Sigma$.
\end{prop}
\begin{pf} D'après la proposition~\ref{exbonnesplaces}, il existe une 
infinité de $v_0$ inertes telles que $q \equiv 1 \pmod{\mu}$ et 
$\psi(T)\equiv q(q^3+1) \mod{\mu^n}$. On en choisit une qui vérifie les
conditions de~\ref{v0}, i.e. qui est 
distincte 
de $v$ et telle que $K_{v_0}$ est hyperspécial. 
On peut alors appliquer le théorème d'augmentation du niveau 
en $v_0$ modulo 
$\mu^n$ de \cite{bg} (théorème 1, page 2). Celui-ci implique l'existence
d'une congruence non triviale modulo $\mu^n$ entre 
$O_{B,J,\pi_\infty,\anneau_\mu}(\psi)$ (l'espace propre de caractère 
propre $\psi$ pour $\Hecke^\Sigma$ sur  $O_{B,J,\pi_\infty,\anneau_\mu}$)
et  $N_{B,J,\pi_\infty,\anneau_\mu}$ donc d'un $f$ non divisible dans 
$N_{B,J,\pi_\infty,\anneau_\mu}$ vérifiant, pour tout $U \in \Hecke^\Sigma$
\begin{eqnarray} \label{cong} U f \equiv \psi(U) f \pmod{\mu^n} \end{eqnarray}

Soit $P$ la projection de 
$N_{B,J,\pi_\infty,\anneau_\mu}$ sur 
$N^+_{B,J,\pi_\infty,\anneau_\mu}$ parallèlement à 
$N^-{B,J,\pi_\infty,\anneau_\mu}$.
La projection $P$ commute à $\Hecke^\Sigma$ car c'est 
un polynôme en $T$. Donc $P(f)$ satisfait aussi la relation~(\ref{cong}).
Par ailleurs, la relation~(\ref{cong}) donne en particulier
$T f \equiv q(q^3+1) f \pmod{\mu^n},$ si bien que $P(f) \equiv f \pmod{\mu^n}$, 
donc $P(f)$ est non divisible. 
Remplaçant $f$ par $P(f)$, la proposition est prouvée.
\end{pf}

Fixons dorénavant $n \geq 1$, et $v_0$ comme dans la proposition précédente.
On pose $\Sigma_0=\Sigma \cup \{v_0\}$ et $\Gamma_E^{\Sigma_0}$ le 
sous-groupe de $\Gamma_E$ correspondant à l'extension maximale de $E$ non
ramifiée hors des places de $E$ divisant les places de $\Sigma_0$. Enfin,
posons 
$\Hecke^{\Sigma_0}=\otimes_{v'\not \in \Sigma_0}\Hecke(G_{v'},K_{v'})$.

Pour $R$ une $\anneau_\mu$-algèbre, notons $\TT_R$ la $R$-algèbre des 
endomorphismes de $N^+_{B,J,\pi_\infty,R}$ engendrée par l'action de 
$\Hecke^{\Sigma_0}$. 
C'est une $R$-algèbre commutative, semi-simple si $R$ est une $L_\mu$-algèbre, 
dont la formation commute à tout changement de base $R \rightarrow R'$ et on a un morphisme 
d'algèbres surjectif $\Hecke^{\Sigma_0} \otimes R \rightarrow \TT_R$. 
La conclusion de la proposition précédente se reformule en 
\begin{prop} \label{reform} 
La réduction modulo $\mu^n$ du caractère $\psi$ de 
$\Hecke^{\Sigma_0}$ se factorise en un caractère, noté $\psi_n$,  
de $\TT_{\anneau_\mu}$ à valeur dans $\anneau_\mu/\mu^n$
\end{prop} 

Par ailleurs, on a
\begin{prop}\label{pseudo}
Il existe un pseudo-caractère $\chi$ de $\Gamma_E^{\Sigma_0}$ à valeur dans 
$\TT_{\anneau_\mu}$ tel que 
\begin{eqnarray} \label{frob} \chi(\frob_{w'}) &=& T_{w'} 
\text{ pour toute place $w'$ divisant une place $v' \not \in \Sigma_0$}\\
\label{iv} \chi(g) &=& \tr r_l((\pi_v)_{E_w})^\ses(g) \ \ \forall g \in I_w
\end{eqnarray}
\end{prop}
\begin{pf}
Il suffit de montrer la proposition pour 
$\TT_{\bar \Q_l}=\TT_{\anneau_\mu} \otimes \bar \Q_l$ puisque $\bar \Q_l$ est
plat sur $\anneau_\mu$. Mais $\TT_{\bar \Q_l}$ est commutative semi-simple et 
de dimension finie, donc est isomorphe à un produit fini de copies 
de $\bar \Q_l$.
Si $\theta$ désigne la projection de $\TT_{\bar \Q_l}$ vers l'un de ces 
facteurs il suffit de prouver l'existence d'un pseudo-caractère $\chi$ de 
$\Gamma_E^{\Sigma_0}$ 
à valeur dans $\bar \Q_l$ vérifiant la propriété~(\ref{iv}) 
et, au lieu de la propriété~(\ref{frob}), la propriété
\begin{eqnarray} \label{frobp} \chi(\frob_{v'}) &=& \theta(T_{v'}) 
\text{ pour toute place $w'$ divisant une place $v' \not \in \Sigma_0$}. \end{eqnarray}

Comme $\bar \Q_l$ est algébriquement clos, 
il existe une forme $0 \neq f \in N^+_{B,J,\pi_\infty,\bar \Q_l}$ propre pour 
$\Hecke^\Sigma$ de caractère propre $\theta$.
D'après~\cite[6.1.2]{bg}, il existe une représentation automorphe 
$\pi'$ attachée à $f$, telle que $\Hecke^{\Sigma_0}$ opère par $\theta$ sur 
$\pi'^{K^{\Sigma_0}}$ (où $K^{\Sigma_0}=\prod_{v' \not \in \Sigma_0} K_{v'}$), et telle 
que $\pi'_{v_0} \simeq \St$ ou $\pi'_{v_0} \simeq \pi^s$.
Soit $\chi$ le caractère de la représentation de $\Gal(\bar E/E)$
sur $\bar \Q_l$ attachée à $\pi'$ ; il vérifie~(\ref{frobp}) 
d'après~\ref{br2}.

Comme $\theta(T) \equiv q(q^3+1) \not \equiv -(q^3+1) \pmod{\mu}$,
on a $\pi_{v_0}=\st$ (cf. \cite[lemme 3.7.3]{bg})
Donc $(\pi'_E)$ est la Steinberg en $v_0$ 
(d'après \cite[prop. 13.2.2(b)]{Roglivre}) qui est de carré
intégrable, et l'on peut appliquer le théorème~\ref{thht} à $\pi'$,
ce qui montre que $\chi(g) = \tr r_l((\pi'_v)_{E_w})^\sss(g)$
pour $g \in D_w$. Mais comme $\Hom_{K_v}(J_v,\pi'_v) \neq 0$, 
le lemme~\ref{type} montre que 
$r_l((\pi'_v)_{E_w})^{\sss}_{|I_w} \simeq r_l((\pi_v)_{E_w})^\sss_{|I_w},$
d'où la formule~(\ref{iv}).
\end{pf}

\subsubsection{Fin de la preuve de la proposition~\ref{propinertie}}

Composons le pseudo-caractère 
$\chi:\ \Gamma_E^{\Sigma_0} \rightarrow \TT_{\anneau_\mu}$ de la 
proposition~\ref{pseudo} avec le caractère $\psi_n$ 
$$\TT_{\anneau_\mu} \rightarrow \anneau_\mu/\mu^n$$ de la 
proposition~\ref{reform}. On obtient un pseudo-caractère 
$$\phi_n=\psi_n \circ \chi:\ 
 \Gamma_E^{\Sigma_0} \rightarrow \anneau_\mu/\mu^n.$$ 
D'après la  formule~(\ref{frob}), et le théorème de Cebotarev,
on a $$\tr \rho_\mu(g) \equiv \phi_n(g) \pmod{\mu^n}\ \ \forall g \in 
\Gamma_E^{\Sigma_0}.$$
Mais d'après~(\ref{iv}), on a 
$$\phi_n(g) \equiv \tr r_l((\pi_v)_{E_w})^\sss(g)\pmod{\mu^n} \ \forall g \in 
I_w.$$
D'où 
$$\tr \rho_\mu(g) \equiv \tr r_l((\pi_v)_{E_w})^\sss(g)\pmod{\mu^n} \ \ \forall g \in I_w $$
Dans cette dernière congruence, les deux membres sont indépendants de $n$ 
(et de $v_0$). Comme celle-ci est valable pour tout $n \geq 1 $, 
on a 
$$\tr \rho_\mu(g) = \tr r_l((\pi_v)_{E_w})^\sss(g) \ \ \forall g \in I_w$$
et la proposition~\ref{propinertie} en 
découle.

\section{Preuve de la proposition~\ref{propiwahori}}

La preuve est très proche de celle de la proposition~\ref{propinertie},
mais un peu plus simple.
Nous indiquons seulement quels sont les changements à faire.

\subsection{Choix d'un niveau et d'un type} 

On prend pour $K_v$ un sous-groupe d'Iwahori $B_v$ de $G_v$
On prend pour $K^v$ un sous-groupe compact ouvert de $\prod_{v'\neq v} G_{v'}$ tel que 
$\pi^{K^v} \neq 0$. On prend $K=B_v K^v$ et $J=1$.

\subsection{Formes anciennes et nouvelles}.

Rien ne change sauf en~\ref{actionhecke} : on définit $\Sigma$ comme l'ensemble
des places $v'$ {\it différentes de $v$} telles que $K_{v'}$ n'est pas
hyperspécial. On définit  $\Hecke^\Sigma$ comme 
$$\left(\bigotimes_{v'\not \in \Sigma,\  v'\neq v} \Hecke(G_{v'},K_{v'}) 
\right) \otimes \ZZ(G_v,B_v),$$ 
et cet anneau agit encore sur les espaces de formes 
$S_{K,J,\pi_\infty,R}$, $O_{B,J,\pi_\infty,R}$, $N^\bullet_{B,J,\pi_\infty,R}$ où $\bullet \in \{\emptyset, +,-\}$. Le lemme~\ref{Nplus} et sa preuve restent
valables sans changement. 

\subsection{Augmentation du niveau}

\subsubsection{}

La proposition~\ref{reform} reste valable sans changement, avec la même preuve.
La proposition~\ref{pseudo} a pour analogue
\begin{prop}\label{pseudo2}
Il existe un pseudo-caractère $\chi$ de $\Gamma_E^{\Sigma_0}$ à valeur dans 
$\TT_{\anneau_\mu}$ tel que 
\begin{eqnarray} \label{frob} \chi(\frob_{w'}) &=& T_{w'}\text{ pour toute place $w'$ divisant une place $v' \not \in \Sigma_0$}
\end{eqnarray}
\end{prop}
\begin{pf}
Comme dans la preuve de la proposition~\ref{reform},
il suffit de montrer, pour tout caractère $\theta$ de $\TT_{\bar \Q_l}$,
 l'existence d'un pseudo-caractère $\chi$ de 
$\Gamma_E^{\Sigma_0}$ à valeur dans $\bar \Q_l$ vérifiant
\begin{eqnarray} \label{frobp2} \chi(\frob_{w'}) &=& \theta(T_{w'}) 
\text{ pour toute place $w'$ divisant une place $v' \not \in \Sigma_0$}.
 \end{eqnarray}

Raisonnant encore comme dans la preuve de la proposition~\ref{reform},
on voit qu'il existe une  une représentation automorphe 
$\pi'$ attachée à $f$, telle que $\Hecke^{\Sigma_0}$ opère par $\theta$ sur 
$\pi'^{K^{\Sigma_0}}$ (où $K^{\Sigma_0}=\prod_{v' \not \in {\Sigma_0}} K_{v'}$), et telle 
que $\pi'_E$ est la Steinberg en la place au-dessus de $v_0$. 

Soit $\chi$ le caractère de la représentation de $\Gal(\bar E/E)$ 
sur $\bar \Q_l$ attachée à $\pi'$ ; 
d'après~\ref{br2}, 
il est non ramifiée hors 
$\Sigma_0 \cup \{v\}$ et l'on a 
$$\chi(\frob_{w'}) = \theta(T_{w'}) 
\text{ pour toute place $w'$ divisant une place $v' \not \in \Sigma$, $v'\neq v$}.$$

D'après le théorème~\ref{thht} appliqué à $\pi'$, 
$\chi(g) = \tr r_l((\pi'_v)_{E_w})^\sss(g)$
pour $g \in D_w$. Donc $\chi$ est aussi non ramifiée en $e$, d'après le point
3) du lemme~\ref{bernstein} et l'on a $$\chi(\frob_w)=\theta(T_w).$$ (La même
chose vaut pour l'autre place $\bar w$ de $E$ au-dessus de $v$ dans le cas où
$v$ est décomposée.)
\end{pf}

\subsubsection{Fin de la preuve de la proposition~\ref{propiwahori}}

Composons le pseudo-caractère 
$\chi:\ \Gamma_E^{\Sigma_0} \rightarrow \TT_{\anneau_\mu}$ de la 
proposition~\ref{pseudo2} avec le caractère $\psi_n$ 
$$\TT_{\anneau_\mu} \rightarrow \anneau_\mu/\mu^n$$ de la 
proposition~\ref{reform}. On obtient un pseudo-caractère 
$$\phi_n=\psi_n \circ \chi:\ 
 \Gamma_E^{\Sigma_0} \rightarrow \anneau_\mu/\mu^n.$$ 
D'après la  formule~(\ref{frob}), et le théorème de Cebotarev,
on a $$\tr \rho_\mu(g) \equiv \phi_n(g) \pmod{\mu^n}\ \ \forall g \in \Gamma_E^{\Sigma_0},$$
et aussi
$$\tr \rho_\mu(g) \equiv \phi_n(g) \pmod{\mu^n}$$
Mais d'après~(\ref{iv}), on a 
$$\tr \rho_\mu(\frob_w) \equiv \psi_n(T_{w}) 
\equiv \tr r_l((\pi_v)_{E_w})^\sss(\frob_w)\pmod{\mu^n} \ \forall g \in I_v$$
la deuxième congruence découlant du lemme~\ref{bernstein}.
Dans cette congruence, les deux membres extrêmes
sont indépendants de $n$ (et de $v_0$). 
Comme elle est valable pour tout $n \geq 1 $, 
on a $$\tr \rho_\mu(\Frob_w) = \tr r_l((\pi_v)_{E_w})^\sss(\frob_w).$$
Ceci est aussi valable pour l'autre place $\bar w$ divisant $w$ (si $v$ est décomposée).  
Comme les représentations $r_l((\pi_v)_{E_w})^\sss$ et $(\rho_\mu)^\sss_{D_w}$
sont non ramifiées (lemme~\ref{bernstein}, point 3),
et ont mêmes caractères centraux (par Cebotarev), elles sont isomorphes, CQFD.

\section{Preuve du théorème principal}

\subsection{Changement de base local}

\subsubsection{}
On note $r$ la représentation $r_l^\ses((\pi_v)_{|E_w})$ du groupe 
$\Gal(\bar \Q_l/ E_w)$.
 Soit $$M=\bar \Q_l^{\ker r}.$$
On a $\Gal(M/E_w)=r(D_w)$.
On dispose de la suite exacte 
$$0 \rightarrow I_w \rightarrow D_w \rightarrow 
\Gal(E_{w}^{nr}/E_w) \rightarrow 0,$$
d'où une suite exacte 
$$0 \rightarrow r(I_w) \rightarrow r(D_w) \rightarrow r(\hat{\Z}) 
\rightarrow 0 $$
Notons que $r(I_w)$ est fini.

\subsubsection{}
Soit $u$ un élément de $D_w = \Gal(\bar \Q_l/ E_w)$.
On suppose que $u \not \in I_w$. On note 
$\langle r(u) \rangle$ l'adhérence du sous-groupe engendré 
par $r(u)$ : sous notre hypothèse, c'est un sous-groupe d'indice fini de 
$r(D_w)$. 

On pose $$N=M^{\langle r(u) \rangle}.$$
On a $\Gal(M/N)=\langle r(u) \rangle$. 
L'extension $N/E_w$ est finie mais n'est pas en général 
galoisienne.

\begin{lemme} \label{mm0}
L'extension $M/N$ est non ramifiée, et son Frobenius est $r(u)$.
\end{lemme}
\begin{pf} Il suffit de montrer que pour toute extension finie galoisienne $M'$
de $N$ contenue dans $M$, on a $M'/N$ non ramifiée, de Frobenius l'image de 
$r(u)$ dans $\Gal(M'/N)$. et c'est clair.
\end{pf}

\subsubsection{} Si $B/A$ est une extension de
corps $p$-adique, on définit le changement de base local de $A$ à $B$ d'une
représentation $\pi$ irréductible lisse de $\Gl_n(A)$, qu'on note $\pi_B$,
par la relation
$$\pi_B := \rec_A^{-1} \circ \res_{A/B} \circ \rec_A(\pi),$$
où $\res_{A/B}$ est la restriction d'une représentation de $W_B$ à $W_A$.

Si la représentation $\pi$ est-elle même le changement de base d'une 
représentation d'un groupe unitaire, par exemple $\pi=(\pi_v)_{E_w}$ 
on notera $(\pi_v)_B$ son changement de base à $B$, au lieu de 
$((\pi_v)_{E_w})_B$, ce qui ne crée pas d'ambiguïté.
\begin{prop} \label{cblocal}
Le changement de base local $(\pi_v)_{N}$ de $(\pi_v)_{E_w}$ à $N$
a une droite invariante par un Iwahori de $\Gl_3(N)$ . 
\end{prop}
\begin{pf} Par définition, $r_l((\pi_v)_{N})$ est la restriction de 
$r_l((\pi_v)_{E_w})$ à $\Gal(\bar \Q_l/N)$. D'après le lemme~\ref{mm0}, cette 
représentation est non ramifiée. La proposition en découle.
\end{pf}
 
Soit $N_0$ la clôture galoisienne de $N$, qui est finie
et résoluble sur $E_w$.

\subsection{Changement de base global}

\begin{lemme}
Il existe un corps de nombres $F'$ vérifiant
\begin{itemize}
\item[i.] $F'$ est totalement réel, 
\item[ii.] $F'/F$ est résoluble et le changement de base $\pi_{EF'}$ de $\Pi_E$ à $EF'$
est cuspidal.
\item[iii] $EF'$ admet une place $w'$ divisant $w$ telle que 
$(EF)_{w'}/E_{w}$ soit isomorphe à $N_0/E_v$.
\end{itemize}
\par 
\end{lemme}
\begin{pf}
Cela résulte aisément de~\cite[théorème 5, page 103]{at}.
Voir~\cite[preuve de la proposition 3.2]{bc} pour plus de détails.
\end{pf}
 
On a $\Gal(N_0/E_w) \subset \Gal(EF'/E)=\Gal(F'/F)$. 
Le sous-groupe $\Gal(N_0/N)$ de $\Gal(N_0/E_w)$ 
s'identifie à un sous-groupe de $\Gal(F'/F)$ et définie donc un sous-corps $H$
de $F'$, contenant $F$, et tel que $(EH)_{w'}/E_w$ est isomorphe à $N/E_w$

Le lemme suivant est démontré dans \cite{harris}.
\begin{lemme}[Harris] Il existe une représentation automorphe cuspidale 
$\pi_{EH}$ de $\Gl_n(EF)$ qui est le changement de base fort de $\pi_E$,
i.e. dont la composante locale $(\pi_{EH})_x$ en toute place $x$ de 
$EH$ divisant une place $y$ de $E$
est le changement de base local de $(\pi_E)_y$ .
\end{lemme}

Par Cebotarev la représentation galoisienne associé à $\pi_{EH}$ est 
$$(\rho_\mu)_{|\Gal(\bar \Q/EH)}.$$

D'après la proposition~\ref{propiwahori}, appliquée à $\pi_{H}$ en la place 
$w'$ de $EH$, on a
$$(\rho_\mu)_{|D_{w'}}^\ses \simeq r_l(\pi_v)_{N} = 
(r(\pi_v)_{E_w})_{|D_{w'}}.$$
Comme $u \in D_{w'}=\Gal(\bar \Q_l/N)$ par définition de $N$, on a en 
particulier
$$ \tr(\rho_\mu)(u) = \tr r(\pi_v)_{E_w}(u) $$
Cette égalité est valable pour tout $u \in D_w-I_w$. Mais par ailleurs,
d'après la proposition~\ref{propinertie}, elle est aussi vrai pour $u \in I_w$.
On a donc montré $\tr (\rho_\mu)_{|D_w} = \tr r_l(\pi_v)_{E_w}$, et le 
théorème. 

\bigskip
{\tiny Joël Bellaïche, jbellaic@math.unice.fr ou jbellaic@ias.edu}

\end{document}